\documentclass[11pt]{article}
\pdfoutput=1

\usepackage[backend=biber,
            bibstyle=phys,
            sorting=none,
            url=false,
            doi=false,
            eprint=true,
            sortcites=true,
            citestyle=numeric-comp]{biblatex}

\usepackage[paperwidth=8.5in, paperheight=11in, portrait, top=1in, bottom=1.in,
 left=1.15in, right=1.15in]{geometry}
\usepackage[utf8]{inputenc}
\usepackage{amssymb,amsthm,amsmath,amsfonts,mathtools,bm}
\usepackage{authblk} 
\usepackage{graphicx,mathrsfs,latexsym,xcolor,subfigure,epsfig}
\usepackage{units}
\usepackage[breaklinks]{hyperref}

\def\sl{{\mathfrak{sl}}}

\def\e{{\textsf{e}}}
\def\P{{{P}}}

\def\p{{\textit{\textsf{P}}}}
\def\f{{\textit{\textsf{F}}}}
\newcommand{\M}[1]{{M_{#1}}}

\numberwithin{equation}{section}
\theoremstyle{plain}
\newtheorem{thm}{Theorem}[section]

\newtheorem{prop}[thm]{Proposition}
\newtheorem{coro}[thm]{Corollary}

\newtheorem{defi}[thm]{Definition}
\newtheorem{rmk}[thm]{Remark}
\newtheorem{example}[thm]{Example}

\title{\bf Racah algebras, the centralizer $Z_n(\sl_2)$\\ and its Hilbert--Poincar\'e series}
\renewcommand*{\Affilfont}{\normalsize\small}
\author[1]{Nicolas Crampé\,}
\author[2]{Julien Gaboriaud\,}
\author[3]{Loïc Poulain d'Andecy\,}
\author[4]{Luc Vinet\,\vspace{.5em}}
\affil[1]{Institut Denis-Poisson CNRS/UMR 7013 - Université de Tours - Université
d'Orléans,
\newline\vspace{.9em}
Parc de Grandmont, 37200 Tours, France.\hfill}
\affil[2,4]{Centre de Recherches Math\'ematiques, Universit\'e de Montr\'eal,
\newline\vspace{.9em}
P.O. Box 6128, Centre-ville Station, Montr\'eal (Qu\'ebec), H3C 3J7, Canada.}
\affil[3]{Laboratoire de mathématiques de Reims UMR 9008,
Université de Reims Champagne-Ardenne,
\newline\vspace{1.5em}
Moulin de la Housse BP 1039, 51100 Reims, France.}
{
 \makeatletter
 \renewcommand\AB@affilsepx{: \protect\Affilfont}
 \makeatother
 \affil[ ]{E-mail addresses}
 \makeatletter
 \renewcommand\AB@affilsepx{, \protect\Affilfont}
 \makeatother
 \affil[1]{crampe1977@gmail.com}
 \affil[2]{julien.gaboriaud@umontreal.ca}
 \affil[3]{loic.poulain-dandecy@univ-reims.fr}
 \affil[4]{vinet@crm.umontreal.ca}
}
\date{\today}

\begin{document}
\maketitle

\hrule
\begin{abstract}
The higher rank Racah algebra $R(n)$ introduced in \cite{DeBieGenestetal2017} is recalled.
A quotient of this algebra by central elements, which we call the special Racah algebra
$sR(n)$, is then introduced. Using results from classical invariant theory, this $sR(n)$
algebra is shown to be isomorphic to the centralizer $Z_{n}(\mathfrak{sl}_2)$ of the
diagonal embedding of $U(\mathfrak{sl}_2)$ in $U(\mathfrak{sl}_2)^{\otimes n}$. This leads
to a first and novel presentation of the centralizer $Z_{n}(\mathfrak{sl}_2)$ in terms of
generators and defining relations. An explicit formula of its Hilbert--Poincaré series is
also obtained and studied.
The extension of the results to the study of the special
Askey--Wilson algebra and its higher rank generalizations is discussed.
\\[.5em]
\noindent{\bf Keywords:} Racah algebra, centralizer, $U(\mathfrak{sl}_2)$, classical
invariant theory, first and second fundamental theorems, Hilbert--Poincar\'e series.
\end{abstract}
\hrule

\section{Introduction}\label{sec:intro}

This paper clarifies the connection between the higher rank Racah algebra $R(n)$ and the
centralizer $Z_n(\sl_2)$ of the diagonal embedding of $U(\sl_2)$ in its $n$-fold tensor
product. The central result of the paper is:\\[.5em]
\centerline{\textit{A particular quotient of the Racah algebra $R(n)$, which we provide in
\eqref{eq:Rnx}, is the centralizer.}}\vspace{.5em}
In the case $n=3$, the centralizer $Z_3(\sl_2)$ contains a number of subalgebras of
interest for mathematics and physics. In representations, it is generated by the so-called
``intermediate'' or ``quadratic'' Casimir elements \cite{LehrerZhang2006}.
These elements realize the relations of the Racah algebra \cite{GenestVinetetal2013}. One
may then wonder if the Racah algebra is the centralizer $Z_3(\sl_2)$ or if more relations
are needed in order to offer a full description of the centralizer by generators and
relations.

A similar story repeats itself for the centralizer $Z_n(\sl_2)$.
One observes that the intermediate Casimir elements in the $n$-fold tensor product
of $U(\sl_2)$ realize the relations of the higher rank Racah algebra \cite{DeBieGenestetal2017}, but two
questions
remain: Is the centralizer generated by these intermediate Casimir elements, and
are there additional relations needed in order to fully describe the centralizer?

As will be seen, these questions are quite close to questions that appear in
classical invariant theory. The issue of finding a generating set (of polynomial
functions, for example) is answered by the First Fundamental Theorem and
the problem of obtaining defining relations between these generators is
answered by the Second Fundamental Theorem.

In this regard, the results that we obtain in Corollary \ref{cor-FFT} and Theorem
\ref{thm:sft_nc} can be seen as non-commutative analogues of the First and Second
Fundamental Theorems.

\subsection{On the Racah algebra}\label{sec:intro_racah}
At this point it would be appropriate to provide some background on the Racah
algebra, which arises in numerous areas of mathematics and physics.

A first approach to the Racah algebra is from the theory of orthogonal polynomials.  The
Racah polynomials are a family of bispectral classical orthogonal polynomials which are characterized
by a difference and recurrence operator \cite{KoekoekLeskyetal2010}. These operators obey
the quadratic algebraic relations of the Racah algebra; this is actually how the algebra
was originally introduced \cite{GranovskiiZhedanov1988} and the reason why it inherited
its name.  Thus, the representation theory of the Racah algebra involves the eponym
polynomials.

The problem of decomposing the tensor product of two irreducible representations of $\sl_2$ in
a direct sum of irreducible representations is known as the Clebsch--Gordan problem of $\sl_2$.
When the three-fold tensor product is considered, there exist two natural decompositions. The question of
finding the overlaps between the two associated bases is called the Racah problem
of $\sl_2$. The intermediate Casimir elements labelling the two mentioned
decompositions realize the Racah algebra \cite{GenestVinetetal2013}, and thus the
overlaps are found to be given in terms of Racah polynomials.  There are also realizations
of the Racah algebra in terms of $U(\sl_2)$ \cite{GranZhe,Koornwinder93, BocktingConradHuand2020,
CrampeShaabankabakibo2021}.

Other instances of the Racah algebra appearing in various contexts include algebraic
combinatorics \cite{Terwilliger2001, GaoWangetal2013},
the theory of double affine Hecke algebras \cite{Huang2021},
and as a symmetry algebra of physical models \cite{KalninsKressetal2005, Post2011a,
KalninsMilleretal2013, GenestVinetetal2014a, GenestVinetetal2014}
Its finite-dimensional irreducible modules have also been classified
\cite{HuangBocktingConrad2020}.

As is the case for various rich structures that have a lot of applications, the
generalization of the Racah algebra is something desirable.
A higher rank Racah algebra has been defined in \cite{DeBieGenestetal2017} by looking
at a $3$D superintegrable system that realizes the Racah algebra as its symmetry algebra
and then generalizing this system to $n$ dimensions. The algebraic relations between the
constants of motion were then computed and used to define abstractly the higher rank
Racah algebra $R(n)$. These relations of $R(n)$ are also verified by the
intermediate Casimir elements that label various direct sum decompositions of the $n$-fold
tensor product of $\sl_2$ irreducible representations \cite{Post2015}.
These higher rank Racah algebras have also been observed in physical models
\cite{KalninsMilleretal2011, DeBieIlievetal2019, LatiniMarquetteetal2020},
interpreted in the framework of Howe duality \cite{GaboriaudVinetetal2019},
and the study of their relation to multivariate Racah polynomials has been
initiated \cite{GeronimoIliev2010, DeBieIlievetal2020, DeBievandeVijver2020}.

\subsection{Outline}\label{sec:intro_outline}
The paper is organized as follows. In section \ref{sec:R3}, the abstract $R(3)$
Racah algebra will be introduced and its quotient by a certain central element will be
presented and named the special Racah algebra $sR(3)$.  Then, the higher rank Racah
algebra $R(n)$ will be presented in Section \ref{sec-Racah} and various properties will be
highlighted.  Section \ref{sec:sRacah} will define the higher rank special Racah algebra
$sR(n)$.  This algebra is a quotient of $R(n)$ by a number of central elements which will
be given precisely. We then come to the main results of the paper. It will be proven in Section \ref{sec:centralizer_proof}
that the special Racah algebra is isomorphic to the diagonal centralizer $Z_n(\sl_2)$ of $U(\sl_2)$ in its $n$-fold tensor
product.  The Hilbert--Poincaré
series of the centralizer will be obtained in Section \ref{sec:hpseries} and its rich
combinatorial properties will be examined.  The PBW basis of the centralizer will then be
given for the first few values of $n$.  A conclusion offering some comments about the
$q$-deformation of these results and the connection with the multivariate Racah
polynomials will end the paper.

\section{The Racah algebra (of rank 1) and a special quotient}\label{sec:R3}

We review the definition of the usual Racah algebra (of rank 1 in our terminology), along
with the algebraic properties we need for the following. We then give the definition of
the special Racah algebra of rank 1, to prepare for the generalization for any rank
defined later in the paper.

\begin{defi}\label{def-R3}
The Racah algebra $R(3)$ of rank $1$ is the associative algebra with generators $\p_{11}$,
$\p_{12}$, $\p_{13}$, $\p_{22}$, $\p_{23}$, $\p_{33}$, $\f_{123}$, and with the following
defining relations for indices $i,j,k$ in $\{1,2,3\}$ and all distinct:
\begin{subequations}\label{eq:defining3}
\begin{align}
 \p_{ii}&\quad \text{is central},\label{eq:r31}\\
 [\p_{ij},\p_{jk}]&=2\f_{ijk},\label{eq:r32}\\
 [\p_{jk},\f_{ijk}]&=\p_{ik}(\p_{jk}+\p_{jj})-(\p_{jk}+\p_{kk})\p_{ij}.\label{eq:r34}
\end{align}
\end{subequations}
where $[A,B]=AB-BA$ is the commutator, and $\p_{ij}$ and $\f_{ijk}$ are defined for any
$i,j,k\in\{1,2,3\}$ by the requirements:
\begin{align}\label{}
 \p_{ij}=\p_{ji}\qquad \text{and}\qquad \f_{ijk}=-\f_{jik}=\f_{jki}\qquad
 \text{for any $i,j,k\in\{1,2,3\}$.}
\end{align}
\end{defi}
We say that $\f_{ijk}$ is antisymmetric whereas $\p_{ij}$ is symmetric. Note that in view
of \eqref{eq:r32}, the element $\f_{123}$ can be removed from the set of generators of
$R(n)$ but it is more convenient to keep it.

Relation (\ref{eq:r32}) is given for any $i,j,k$. If $i,j,k$ are not ordered, the symmetry
properties of $\p$ and $\f$ are to be used. For example,
\begin{align}\label{}
 [\p_{23},\p_{13}]=[\p_{23},\p_{31}]=2\f_{231}=2\f_{123},
\end{align}
where we used the symmetry of $\p$, the antisymmetry of $\f$ and relation \eqref{eq:r32}
for $(i,j,k)=(2,3,1)$. Similar comments apply to relation \eqref{eq:r34}, which is
given for any distinct $i,j,k$ and not only the case $(i,j,k)=(1,2,3)$.

Let us introduce the notion of determinant for matrices with non-commuting entries. If $A$
is a $n\times n$ matrix with entries $A_{i,j}$ ($1 \leq i,j \leq n$), we define the
symmetrized determinant of $A$ as follows
\begin{align}\label{eq:sym_det}
 \det A:=\frac{1}{n!} \sum_{\rho, \sigma\in
S_{n}}sgn(\rho) sgn(\sigma)A_{\rho(1),\sigma(1)}A_{\rho(2),\sigma(2)}\dots
 A_{\rho(n),\sigma(n)},
\end{align}
where $S_n$ is the permutation group of $n$ elements and $sgn(\sigma)$ is the signature of
$\sigma$.
For commuting entries, it is the usual definition of the determinant of a matrix.
We define also the following $3\times 3$ matrix
\begin{align}\label{eq:Pijkabc}
 \p_{ijk}^{\,abc}=
 \begin{pmatrix}
  \p_{ia} & \p_{ib} & \p_{ic}\\
  \p_{ja} & \p_{jb} & \p_{jc}\\
  \p_{ka} & \p_{kb} & \p_{kc}
 \end{pmatrix}.
\end{align}
There exist in addition to $P_{ii}$, other central elements in $R(3)$ given in the following proposition:
\begin{prop}\label{pr:R3}
The following elements are central in $R(3)$\begin{align}
 &Q_{3}=\p_{12}+\p_{13}+\p_{23},\\
 &\begin{aligned}
  w_{ijk}:&={\f_{ijk}}^2 + \tfrac12\det( \p_{ijk}^{\,ijk} )\\
  &\qquad
   -\tfrac13\left( \{\p_{ij},\p_{ik}\} + \{\p_{ij},\p_{jk}\} + \{\p_{ik}\p_{jk}\}
                 + \p_{ij}\p_{kk} + \p_{ik}\p_{jj} + \p_{jk}\p_{ii} \right),
  \end{aligned}\label{}
\end{align}
for $1\leq i,j,k \leq 3$ distinct and where the anticommutator is defined as
$\{A,B\}=AB+BA$. Moreover it is observed that $w_{ijk}$ is symmetric \textit{i.e.}
\begin{align}\label{eq:wijkjik}
 w_{ijk}=w_{jik}=w_{jki}.
\end{align}
\end{prop}
\proof
To show that an element is central, it is enough to show that it commutes with
$\p_{12},\p_{13},\p_{23}$ (since the $\p_{ii}$'s are central and
$2\f_{123}=[\p_{12},\p_{23}]$). For $Q_3$, this is an easy verification.

The symmetry of $w_{ijk}$ is immediate. So it remains to show the centrality of $w_{123}$.
By symmetry of the algebra under renaming of the indices, only the commutation with one
element, say $P_{23}$, needs to be checked. This is a direct calculation using the
defining relations of the algebra.
\endproof

\begin{rmk}
It should be stressed that the element $w_{123}$ is essentially the known Casimir
element of the Racah algebra \cite{GranovskiiZhedanov1988}.
\end{rmk}

\paragraph{The special Racah algebra.}
The special Racah algebra $sR(3)$ is defined from the Racah algebra $R(3)$ by fixing the
value of this non-trivial central element.
\begin{defi}
The special Racah algebra $sR(3)$ of rank $1$ is the quotient of $R(3)$ by
\begin{align}\label{eq:R3x}
 w_{123}=0.
\end{align}
\end{defi}
\begin{rmk}
This relation \eqref{eq:R3x} is akin to the relation that expresses the Casimir
element of the Racah algebra in terms of its central elements, see equation (3.4) in
\cite{GenestVinetetal2013} for instance.
\end{rmk}
\begin{rmk}
It follows from \eqref{eq:wijkjik} that $w_{213}$ and other $w_{ijk}$ obtained from
permutations of the indices are null in $sR(3)$.
\end{rmk}
The appellation ``special'' is used in the same way as in \cite{CrampeFrappatetal2021},
where a quotient of the Askey--Wilson algebra by fixing the value of a central element
expressed as a determinant was denoted as the ``special Askey--Wilson algebra'' (this was
inspired by the nomenclature of Lie groups).

\section{The higher rank Racah algebras}\label{sec-Racah}

Following \cite{DeBieGenestetal2017}, we consider the following definition of the Racah
algebra (of any rank). Note that we consider it as an abstract algebra defined by
generators and relations. It will be clear that for $n=3$ we recover Definition
\ref{def-R3} for the rank one Racah algebra $R(3)$.
\begin{defi}
The Racah algebra $R(n)$ of rank $n-2$ is the associative algebra with generators:
\begin{align}\label{}
 \p_{ij},\ 1\leq i\leq j\leq n\qquad \text{and}\qquad\f_{ijk},\ 1\leq i<j<k\leq n,
\end{align}
and the defining relations are, for all possible indices $i,j,k,l,m$ in $\{1,\dots,n\}$:
\begin{subequations}\label{eq:definingn}
\begin{align}\label{}
 \p_{ii}&\quad \text{is central},\label{eq:rn0}\\
 [\p_{ij},\p_{k\ell}]&=0\hspace{6em}
 \text{if both $i,j$ are distinct from $k,\ell$},\label{eq:rn1}\\
 [\p_{ij},\p_{jk}]&=2\f_{ijk},\label{eq:rn2}\\
 [\p_{jk},\f_{ijk}]&=\p_{ik}(\p_{jk}+\p_{jj})-(\p_{jk}+\p_{kk})\p_{ij},\label{eq:rn3}\\
 [\p_{k\ell},\f_{ijk}]&=\p_{ik}\p_{j\ell}-\p_{i\ell}\p_{jk},\label{eq:rn4}\\
 [\f_{ijk},\f_{jk\ell}]&=-(\f_{ij\ell}+\f_{ik\ell})\p_{jk},\label{eq:rn5}\\
 [\f_{ijk},\f_{k\ell m}]&=\f_{i\ell m}\p_{jk}-\f_{j\ell m}\p_{ik},\label{eq:rn6}
\end{align}
\end{subequations}
where in each relation all indices involved are distinct and
$\p_{ij}$ and $\f_{ijk}$ are defined by:
\begin{align}\label{}
 \p_{ij}=\p_{ji}\qquad \text{and}\qquad \f_{ijk}=-\f_{jik}=\f_{jki}\qquad
 \text{for any $i,j,k\in\{1,\dots,n\}$.}
\end{align}
\end{defi}
We say that $\f_{ijk}$ (or simply $\f$) is antisymmetric whereas $\p_{ij}$ is symmetric.
The same comments as for the algebra $R(3)$ after Definition \ref{def-R3} apply here. In
particular, in view of \eqref{eq:rn2}, the elements $\f_{ijk}$ can be removed from the set
of generators of $R(n)$ but it is more convenient to keep them. Note that to check that a
certain element $X$ is central in $R(n)$, it is enough to check that it commutes with all
generators $\p_{ij}$.

\begin{rmk}
The form of the defining relations above is very symmetrical, and this is quite useful in
practice. Namely, for any permutation $\pi$ of $\{1,\dots,n\}$, the corresponding renaming of
the generators $(\p_{ij},\f_{ijk})\mapsto (\p_{\pi(i),\pi(j)},\f_{\pi(i),\pi(j),\pi(k)})$ is
an automorphism of the algebra. So when checking a relation in the algebra $R(n)$, it is
enough to do it for a chosen set of indices. This property will be used in the proofs.
\end{rmk}

Recall that $\p_{ii}$ is central in $R(n)$. It is easy to show using
(\ref{eq:rn1})-(\ref{eq:rn2}) that the following element is central in $R(n)$:
\begin{align}\label{Qn}
 Q_{n}=\sum_{1\leq i<j \leq n} \p_{ij}\ .
\end{align}
We now introduce some elements of $R(n)$ that will later play an important part.
Recall the definition of $\det(\p_{ijk}^{abc})$ \eqref{eq:Pijkabc}
formulated in the preceding section. We define the following elements in $R(n)$:
\begin{align}
&\begin{aligned}
  w_{ijk}&:={\f_{ijk}}^2 + \tfrac12\det( \p_{ijk}^{\,ijk} )\\
  &\qquad
   -\tfrac13\left( \{\p_{ij},\p_{ik}\} + \{\p_{ij},\p_{jk}\} + \{\p_{ik}\p_{jk}\}
                 + \p_{ij}\p_{kk} + \p_{ik}\p_{jj} + \p_{jk}\p_{ii} \right),
 \end{aligned}\label{eq:defw} \\[0.4em]
& x_{ijk\ell}: = \f_{ijk}\f_{jk\ell}+\tfrac12\det(\p_{ijk}^{jk\ell})
 + \tfrac12(\f_{ij\ell}+\f_{ik\ell})\p_{jk}
 - \tfrac13(\p_{ij}\p_{k\ell} + \p_{ik}\p_{j\ell} + \p_{i\ell}\p_{jk}),
 \label{eq:defx}\\[0.4em]
& y_{ijk\ell m}: = \f_{ijk}\f_{k\ell m} + \tfrac12\det(\p_{ijk}^{k\ell m})
               + \tfrac12(\f_{ij\ell}\p_{km} - \f_{ijm}\p_{k\ell}),
 \label{eq:defy}\\[0.4em]
& z_{ijk\ell mp} := \f_{ijk}\f_{\ell mp} + \tfrac12\det(\p_{ijk}^{\ell mp}),
 \label{eq:defz}
\end{align}
where indices $i,j,k,\ell,m,p\in\{1,\dots,n\}$ are all distinct. Only the element
$w_{ijk}$ appears in the Racah algebra $R(3)$ since there are not enough different
indices for the other elements.

\subsection{The Racah algebra $R(4)$}\label{sec:R4}

The algebra $R(3)$ was previously studied in Section \ref{sec:R3} so let us consider now
the case $n=4$. According to the definition given above, there are $10$ generators
$\p_{ij}$ of $R(4)$ and $4$ generators $\f_{ijk}$.  Of the relations
(\ref{eq:rn1})--(\ref{eq:rn5}) only those involving no more than 4 different indices are necessary
here.

We already know that the elements $\p_{ii}$ and $Q_{4}$ given in
\eqref{Qn} are central in $R(4)$. The following proposition gives less immediate consequences of the defining
relations of $R(4)$, and in particular identifies the elements introduced
in (\ref{eq:defw})--(\ref{eq:defx}) as central elements.
\begin{prop}\label{pr:R4}
The following assertions are true in $R(4)$:
\begin{itemize}
\item For $1\leq a \leq 4$, the following relations hold:
 \begin{align}\label{eq:r4d5}
 & \p_{a1}\f_{234}-\p_{a2}\f_{134}+\p_{a3}\f_{124}-\p_{a4}\f_{123}=0.
\end{align}
\item For distinct $i,j,k\in\{1,2,3,4\}$, the elements $w_{ijk}$ are symmetric
($w_{ijk}=w_{jik}=w_{jki}$) and are central in $R(4)$.
\item For distinct $i,j,k,\ell\in\{1,2,3,4\}$, the elements $x_{ijkl}$ are symmetric
($x_{\sigma(i)\sigma(j)\sigma(k)\sigma(\ell)}=x_{ijk\ell}$ for $\sigma\in S_4$) and are
central in $R(4)$.
\end{itemize}
\end{prop}
\proof
All these statements are proven by invoking the associativity of the algebra.
Here is what is meant by that.
Suppose that we have a word $CBA$, for $A$, $B$, $C$ some generators, that we want to
reorder into the form $ABC$.
This is done by using the defining relations of the algebra
\eqref{eq:definingn}. We decide, for example, to bring all $\p$'s to the left of
the $\f$'s, and to order the $\f$'s and the $\p$'s between themselves in the
lexicographical ordering of their indices. There are two ways to proceed: one can first
start by reordering the pair $(CB)$ into $(BC)+\text{some terms}$, or one could instead
start by reordering the pair $(BA)$ into $(AB)+\text{some other terms}$.  We denote
symbolically the difference at the end of these two computations by
\begin{align}\label{}
 \bigl(C(BA)-(CB)A\bigr),
\end{align}
and this must be identically $0$ by the associativity of the algebra.

Let us first prove \eqref{eq:r4d5}.
In the present case, we shall look at the word $CBA=\f_{ij\ell}\p_{k\ell}\p_{i\ell}$,
for $i,j,k,\ell$ all distinct, and compute
\begin{align}\label{}
 \tfrac12\bigl((\f_{ij\ell}\p_{k\ell})\p_{i\ell}-\f_{ij\ell}(\p_{k\ell}\p_{i\ell})\bigr).
\end{align}
Using relations \eqref{eq:definingn}, this can be brought to the form
\begin{align}\label{}
 \p_{ii}\f_{jk\ell}-\p_{ij}\f_{ik\ell}+\p_{ik}\f_{ij\ell}-\p_{i\ell}\f_{ijk}.
\end{align}
By the argument above, this expression has to be zero.
Then, choosing $(i,j,k,\ell)=(1,2,3,4)$, $(2,3,4,1)$, $(3,4,1,2)$ or $(4,1,2,3)$, we
recover \eqref{eq:r4d5} with $a=1,2,3,4$ respectively.

The proof that $w_{ijk}$ is symmetric and that it commutes with all $\p_{ab}$ with
$a,b\in\{i,j,k\}$ was already done in Section \ref{sec:R3} for the algebra $R(3)$, and is
still valid here.
Using the symmetry of $w_{ijk}$ and the symmetry of the algebra under renaming of the
indices, to prove that $w_{ijk}$ is central, it is enough to prove for example that
$[\p_{34},w_{123}]=0$.

This is done by making use of \eqref{eq:r4d5} and reducing the calculation to
\begin{align}\label{}
 [\p_{34},w_{123}]=\bigl((\f_{124}+2\f_{134})\f_{123}\bigr)\p_{23}
 -(\f_{124}+2\f_{134})\bigl(\f_{123}\p_{23}\bigr)
\end{align}
which is identically zero by the associativity of the algebra.

For the symmetry of $x_{ijk\ell}$, the particular case of $x_{ijk\ell}=x_{ikj\ell}$ is
immediate using the symmetries of $\p$ and $\f$. To complete the proof of the symmetry
properties of $x$, it remains to show that $x_{jk\ell i}=x_{ijk\ell}$. Using the symmetry
of the algebra $R(4)$ under renaming of the indices, it is enough to check that for
example $x_{2341}=x_{1234}$. Substituting from the definition of $x_{ijk\ell}$, one has
\begin{align}\label{eq:diffxx}
\begin{aligned}
 x_{2341}-x_{1234}&=(\f_{134}-\f_{123})\f_{234}
  + \tfrac{1}{2}\bigl(\det(\p_{134}^{234})-\det(\p_{123}^{234})\bigr)\\
  &- \tfrac{1}{2}\bigl(\p_{34}(\f_{123}+\f_{124})+\p_{23}(\f_{134}+\f_{124})\bigr).
\end{aligned}
\end{align}
Now, looking at $CBA=\f_{234}\f_{124}\p_{34}$ and making use of \eqref{eq:r4d5}, one
computes
\begin{align}\label{eq:diffxx2}
\begin{aligned}
 \tfrac{1}{2}\bigl(\f_{234}(\f_{124}\p_{34})-(\f_{234}\f_{124})\p_{34}\bigr)
 &=(\f_{134}-\f_{123})\f_{234}
  + \tfrac{1}{2}\bigl(\det(\p_{134}^{234})-\det(\p_{123}^{234})\bigr)\\
  &\quad- \tfrac{1}{2}\bigl(\p_{34}(\f_{123}+\f_{124})+\p_{23}(\f_{134}+\f_{124})\bigr)\ .
\end{aligned}
\end{align}
By the associativity of the algebra (see above), this expression has to be zero. This
completes the proof of the symmetry of $x_{ijk\ell}$ since the right hand sides of
(\ref{eq:diffxx}) and (\ref{eq:diffxx2}) are the same.

Using the symmetry of $x_{ijk\ell}$ and the symmetry of the algebra under renaming of the
indices, the proof that $x_{ijk\ell}$ is central reduces to proving that for example
$[x_{1234},\p_{23}]=0$ which is also done by a direct computation using expression \eqref{eq-xexp} of $x_{1234}$.
\endproof

\begin{rmk}
 The elements $w_{ijk}$ and $x_{ijk\ell}$ can be equivalently given by the following formulae
\begin{align}
 w_{123}&={\f_{123}}^2 -\f_{123}\p_{13}-\p_{12}(\p_{13}+\p_{23}+\p_{33}) +\frac{1}{2} \sum_{\sigma\in S_3}sgn(\sigma) \p_{\sigma(1)1} \p_{\sigma(2)2} \p_{\sigma(3)3},\label{eq-wexp}
\end{align}
and
\begin{align}
 x_{1234}&=  \f_{123}\f_{234} -\f_{123}\p_{24}+\f_{124}\p_{23}+\f_{134}\p_{23} -\p_{14}\p_{23}
 +\frac{1}{2} \sum_{\sigma\in S_3}sgn(\sigma) \p_{\sigma(1)2} \p_{\sigma(2)3} \p_{\sigma(3)4}. \label{eq-xexp}
\end{align}

\end{rmk}

\subsection{The Racah algebra $R(n)$ for any $n$}

Let now $n$ be any positive integer. We already know that we have in $R(n)$ central
elements $\p_{ii}$ and $Q_{n}$ given in (\ref{Qn}). Building upon all that has been proven
up to now, we have the following final proposition about the Racah algebra $R(n)$.
\begin{prop}\label{pr:R6}
The following assertions are true in $R(n)$:
\begin{itemize}
\item The relations below hold for $1\leq a \leq n$ and $1\leq
i< j< k < \ell \leq n$:
\begin{align}\label{eq:rnd5}
 & \p_{ai} \f_{jk\ell}-\p_{aj} \f_{ik\ell}+\p_{ak} \f_{ij\ell}-\p_{a\ell} \f_{ijk}=0.
\end{align}
\item For distinct $i,j,k\in\{1,\dots,n\}$, the elements $w_{ijk}$ are symmetric and
central in $R(n)$.
\item For distinct $i,j,k,\ell\in\{1,\dots,n\}$, the elements $x_{ijk\ell}$ are symmetric
and central in $R(n)$.
\item For distinct $i,j,k,\ell,m\in\{1,\dots,n\}$, the elements $y_{ijk\ell m}$ are null
in $R(n)$:
\begin{align}\label{}
 y_{ijk\ell m}=0\qquad \text{for all distinct $i,j,k,\ell,m$.}
\end{align}
\item For distinct $i,j,k,\ell,m,p\in\{1,\dots,n\}$, the elements $z_{ijk\ell mp}$ are null in $R(n)$:
\begin{align}\label{}
 z_{ijk\ell mp}=0\qquad \text{for all distinct $i,j,k,\ell,m,p$.}
\end{align}
\end{itemize}
\end{prop}
\proof
If $a$ is equal to $i,j,k$ or $\ell$, relation \eqref{eq:rnd5} only involves 4 indices,
and so its validity follows directly from the Proposition \ref{pr:R4} concerning the
algebra $R(4)$. To prove the case when $a$ is different of $i,j,k$ and $\ell$, first
compare \eqref{eq:rn6} for $(i,j,k,\ell,m)$ equals to $(1,2,5,3,4)$ and $(3,4,5,1,2)$.
This leads to the identity
\begin{align}\label{eq:r5d5s}
 & \p_{15} \f_{234}-\p_{25} \f_{134}+\p_{35} \f_{124}-\p_{45} \f_{123}=0.
\end{align}
The symmetry of the algebra under renaming of the indices suffices to
complete the proof of \eqref{eq:rnd5}.

Regarding the statements about $w_{ijk}$ and building on the verifications made in $R(3)$
and $R(4)$, it remains only to check that $w_{ijk}$ commutes with $\p_{\ell m}$ when
$\ell,m\notin\{i,j,k\}$. This is immediate since from the defining relations, any two
elements ($\p$'s or $\f$'s) with no index in common commute.

Concerning $x_{ijk\ell}$, building on the proof of the preceding subsection, it remains to
show, for example, that
\begin{align}
 [x_{1234}, \p_{45}]=0.
\end{align}
This is shown using $y_{ijk\ell m}=0$ (which is proven below) as well as \eqref{eq:rnd5}.
For $x_{ijk\ell}$, it remains only to check that it commutes with $\p_{mp}$ when
$m,p\notin\{i,j,k,\ell\}$, and this is immediate.

To prove that $y_{ijk\ell m}$ is zero, we first look at the case
$(i,j,k,\ell,m)=(1,2,3,4,5)$. Write
\begin{align}\label{}
 \tfrac{1}{2}\bigl(\f_{245}(\f_{123}\p_{23})-(\f_{245}\f_{123})\p_{23}\bigr).
\end{align}
It is seen that this is equal to $y_{12345}$, by making use of \eqref{eq:rnd5}.  Invoking
associativity, it follows that $y_{12345}=0$. Since this can be repeated for all other
combinations of distinct indices $i,j,k,\ell,m$ it is done.

For the nullity of $z_{ijk\ell mp}$, we invoke again the associativity of the algebra.
Looking at
\begin{align}\label{}
 \tfrac{1}{2}\bigl((\f_{356}\f_{234})\p_{12}-\f_{356}(\f_{234}\p_{12})\bigr)
\end{align}
it is seen that this is equal to $z_{123456}$. Thus it follows that
\begin{align}\label{}
 z_{123456}=0.
\end{align}
A similar reasoning can be repeated for different indices to complete the proof.
\endproof

\section{The special Racah algebra $sR(n)$ }\label{sec:sRacah}

After the preliminary discussion of the Racah algebra $R(n)$, we are now ready to define the
special Racah algebra for any $n$. This is a generalization to arbitrary rank of the
special Racah algebra $sR(3)$ from Section \ref{sec:R3}.
\begin{defi}
The special Racah algebra $sR(n)$ of rank $n-2$ is the quotient of $R(n)$ by all
\begin{align}\label{eq:Rnx}
 w_{ijk}=0,\qquad
 x_{abcd}=0,
\end{align}
such that $1\leq i<j<k\leq n$ and $1\leq a<b<c<d\leq n$.
\end{defi}
Since $sR(n)$ is the algebra involved for the study of the centralizer in the next
section, we collect here the generators and defining relations, to give an explicit
definition without reference to the Racah algebra (we keep the same name for the
generators, which is justified since this is a quotient, this should not lead to any
ambiguity).

\begin{defi}[Equivalent definition] The special Racah algebra $sR(n)$, of rank $n-2$ is
the associative algebra with generators:
\begin{align}\label{}
 \p_{ij},\ 1\leq i\leq j\leq n\qquad \text{and}\qquad\f_{ijk},\ 1\leq i<j<k\leq n.
\end{align}
To give the defining relations, first we define $\p_{ij}$ and $\f_{ijk}$ by:
\begin{align}\label{}
 \p_{ij}=\p_{ji}\qquad \text{and}\qquad \f_{ijk}=-\f_{jik}=\f_{jki}\qquad
 \text{for any $i,j,k\in\{1,\dots,n\}$.}
\end{align}
The defining relations are, for all possible indices $i,j,k,\ell,m$ in $\{1,\dots,n\}$:
\begin{subequations}\label{eq:definingsrn}
\begin{align}\label{}
 \p_{ii}&\quad \text{is central},\label{eq:srn0}\\
 [\p_{ij},\p_{k\ell}]&=0\hspace{6em}\text{if both $i,j$ are distinct from $k,\ell$},
 \label{eq:srn1}\\
 [\p_{ij},\p_{jk}]&=2\f_{ijk},\label{eq:srn2}\\
 [\p_{jk},\f_{ijk}]&=\p_{ik}(\p_{jk}+\p_{jj})-(\p_{jk}+\p_{kk})\p_{ij},\label{eq:srn3}\\
 [\p_{k\ell},\f_{ijk}]&=\p_{ik}\p_{j\ell}-\p_{i\ell}\p_{jk},\label{eq:srn4}\\
 [\f_{ijk},\f_{jk\ell}]
 &=\f_{jk\ell}\p_{ij}-\f_{ijk}\p_{j\ell}-\f_{ik\ell}(\p_{jk}+\p_{jj}),\label{eq:srn5}\\
 [\f_{ijk},\f_{k\ell m}]&=\f_{i\ell m}\p_{jk}-\f_{j\ell m}\p_{ik},\label{eq:srn6}
\end{align}
along with
\begin{align}
 {\f_{ijk}}^2 + \tfrac12\det( \p_{ijk}^{ijk} )&=
 \tfrac13\left( \{\p_{ij},\p_{ik}\} + \{\p_{ij},\p_{jk}\} + \{\p_{ik}\p_{jk}\}
               + \p_{ij}\p_{kk} + \p_{ik}\p_{jj} + \p_{jk}\p_{ii} \right),
 \label{eq:srn7}
\end{align}
for $1\leq i<j<k\leq n$, and
\begin{align}
 \f_{ijk}\f_{jk\ell}+\tfrac12\det(\p_{ijk}^{jk\ell})&=
 - \tfrac12(\f_{ij\ell}+\f_{ik\ell})\p_{jk}
 + \tfrac13(\p_{ij}\p_{k\ell} + \p_{ik}\p_{j\ell} + \p_{i\ell}\p_{jk}),\label{eq:srn8}
\end{align}
for $1\leq i<j<k<\ell\leq n$.
\end{subequations}
\end{defi}

\paragraph{Consequences of the relations.}
From the results of the preceding sections, we know that relations
(\ref{eq:srn7})-(\ref{eq:srn8}) for any distinct $i,j,k,\ell$ are automatically verified,
along with the relations:
\begin{subequations}\label{eq:consequences}
\begin{align}
 \f_{ijk}\f_{k\ell m}+\tfrac12\det(\p_{ijk}^{k\ell m})&=
 \tfrac12(\f_{ijm}\p_{k\ell} - \f_{ij\ell}\p_{km})\\
 \f_{ijk}\f_{\ell mr}+\tfrac12 \det(\p_{ijk}^{\ell mr})&=0,\\
 \p_{ai} \f_{jk\ell}-\p_{aj}\f_{ik\ell}+\p_{ak} \f_{ij\ell}-\p_{a\ell} \f_{ijk}& =0.
 \label{eq:srn11}
\end{align}
\end{subequations}
for distinct $i,j,k,\ell,m,r\in\{1,\dots,n\}$ and for any $a\in\{1,\dots,n\}$. These
relations, although satisfied, do not have to be included in the set of defining
relations.

\begin{example}
The special Racah algebra $sR(4)$ of rank $2$ is the quotient of $R(4)$ by
\begin{align}\label{eq:R4x}
 w_{123}=0,\quad w_{124}=0,\quad w_{134}=0,\quad w_{234}=0\quad
 \text{and}\quad x_{1234}=0.
\end{align}
\end{example}
\begin{rmk}
Some analogues of the relations of $sR(4)$ (excluding the ones of the type
\eqref{eq:srn11})
were obtained in a particular realization in \cite{KalninsMilleretal2011}.
For the higher rank case of $sR(n)$, analogous relations were also observed in
a certain realization in \cite{LatiniMarquetteetal2020}, once again excluding the ones of
the type \eqref{eq:srn11}.
\end{rmk}

\section{Isomorphism between the centralizer $Z_n(\sl_2)$ and the special Racah algebra $sR(n)$}
\label{sec:centralizer_proof}

The goal of this section is to connect the (higher rank) special Racah algebra
introduced and characterized in the previous two sections with the centralizer
$Z_n(\sl_2)$ of the diagonal action of $U(\sl_2)$ in $U(\sl_2)^{\otimes n}$.  This will
provide an a posteriori justification for the quotient that was chosen to go from the
Racah algebra $R(n)$ to the special Racah algebra $sR(n)$: as will be shown, this quotient
is precisely the one that leads to an algebra isomorphic to the centralizer
$Z_n(\sl_2)$.

\subsection{Centralizer $Z_n(\sl_2)$ of the diagonal action $U(\sl_2)$ into
$U(\sl_2)^{\otimes n}$ and the algebra of polarized traces}

We here define the centralizer associated to the Lie algebra $\sl_2$.
The generators of $\sl_2$ are $e_{ij}$, $i,j\in\{1,2\}$ obeying the defining relations
\begin{align}\label{}
 [e_{ij},e_{k\ell}]=\delta_{jk}e_{i\ell}-\delta_{\ell i}e_{kj},\qquad
 e_{11}+e_{22}=0.
\end{align}
We denote by $U(\sl_2)$ the universal enveloping algebra of $\sl_2$. Its Casimir element
is given by
\begin{align}\label{}
 C={e_{11}}^{2}-e_{11}+e_{12}e_{21}.
\end{align}
Now consider the tensor product of $n$ copies of $U(\sl_2)$ and define the following
notation for its generators
\begin{align}\label{}
 e_{ij}^{(a)}=1^{\otimes(a-1)}\otimes e_{ij}\otimes1^{\otimes(n-a)}.
\end{align}
The diagonal embedding of $U(\sl_2)$ in its $n$-fold tensor product is given by
\begin{align}\label{}
\begin{aligned}
 \delta:~U(\sl_2)&\to U(\sl_2)^{\otimes n}\\
 ~e_{ij}&\mapsto\sum_{a=1}^{n}e_{ij}^{(a)}.
\end{aligned}
\end{align}
There is a natural degree-preserving action of $\sl_2$ on $U(\sl_2)^{\otimes n}$ given
by composing the diagonal embedding $\delta$ followed by the adjoint action. On the
generators, it is given by
\begin{align}\label{eq:actionUsl2n}
 e_{ij}\cdot e_{k\ell}^{(a)}=\delta_{jk}e_{i\ell}^{(a)}-\delta_{\ell i}e_{kj}^{(a)}.
\end{align}
We then define the centralizer $Z_n(\sl_2)$ of the diagonal embedding of $U(\sl_2)$ in
$U(\sl_2)^{\otimes n}$ as the kernel of this $\sl_2$ action
\begin{align}\label{eq:centralizerCn}
 Z_n(\sl_2)=\left\{X\in U(\sl_2)^{\otimes n}~|~g\cdot X=[\delta(g),X]=0 \quad
 \forall g\in U (sl_2)\right\}
\end{align}
or in other words, as the set of elements in $U(\sl_2)^{\otimes n}$ that commute with
the diagonal embedding of $U(\sl_2)$.

Let us also define the polarized traces (the summation convention is assumed):
\begin{align}\label{eq:polarizedtraces}
 T^{(a_1,\dots,a_d)}=e_{i_2i_1}^{(a_1)}e_{i_3i_2}^{(a_2)}\dots e_{i_1i_d}^{(a_d)},\qquad
 a_1,\dots,a_d\in\{1,\dots,n\}.
\end{align}
It is seen by a direct computation that these elements are in the centralizer
$Z_n(\sl_2)$.
\begin{rmk}
It is easily checked from the definition \eqref{eq:polarizedtraces} that
$T^{(a_{1},a_{2})}=T^{(a_{2},a_{1})}$ and that $T^{(a_{1},a_{2},a_{3})}$ is antisymmetric
in its indices $a_{1},a_{2},a_{3}$, i.e.
$T^{(a_{1},a_{2},a_{3})}=T^{(a_{2},a_{3},a_{1})}=-T^{(a_{2},a_{1},a_{3})}$.
\end{rmk}
\begin{rmk}
A number of papers in the literature \cite{GenestVinetetal2013, Post2015}
realize the Racah algebra with the so-called
``intermediate Casimir'' elements $C_{i}$, $C_{ij}$. These elements are given by
\begin{align}\label{}
 C_{i}=1^{\otimes(i-1)}\otimes C\otimes1^{\otimes(n-i)},\qquad
 C_{ij}=1^{\otimes(i-1)}\otimes C_{(1)}\otimes1^{\otimes(j-i-1)}\otimes
                                C_{(2)}\otimes1^{\otimes(n-j)},
\end{align}
where $\Delta(e_{ij})=e_{ij}\otimes1+1\otimes e_{ij}$ and we denote
$\Delta(C)=C_{(1)}\otimes C_{(2)}$ in Sweedler's notation.
Then the $T^{(i,i)}$ and $T^{(i,j)}$ can be expressed in terms of these intermediate
Casimir elements as follows
\begin{align}\label{}
 T^{(i,i)}=2C_{i},\qquad T^{(i,j)}=C_{ij}-C_{i}-C_{j}.
\end{align}
\end{rmk}

\subsection{Elements of classical invariant theory \label{sec-inv}}

We now present results from classical invariant theory about the algebra of polynomial
functions on
\begin{align}\label{}
 \underbrace{\sl_2\times\dots\times\sl_2}_{\text{$n$ factors}}\equiv\sl_2^n,
\end{align}
that are invariant under simultaneous conjugations by $SL(2)$.
For $G$ elements of $SL(2)$, these actions on a polynomial
function of $\sl_2^n$ are given by:
\begin{align}\label{eq:simultaneous_conjug}
 G\cdot f(\M{1},\dots,\M{n})
 =f(G^{-1}\M{1}G,\dots,G^{-1}\M{n}G),
\end{align}
for $\M{i}\in\sl_2$.
The first fundamental theorem of classical invariant theory states that:
\begin{thm}[see \cite{Procesi1976,Razmyslov1974,Sibirskii1968} or \cite{Drensky2006} and
references therein]\label{th:fft}
The algebra $\mathbb{C}[\sl_2^n]^{inv}$ of polynomial functions on $\sl_2^n$ that are
invariant under simultaneous conjugations by $SL(2)$ elements is generated by the functions
\begin{align}\label{eq:generatorsinvtfct}
\mathfrak{T}^{(a_1,\dots, a_d)}\ :\  (\M{1},\dots,\M{n})\mapsto Tr(\M{a_1}\dots \M{a_d})
\end{align}
for $\M{i}\in\sl_2$, $d\geq2$ and $a_1,\dots,a_d\in\{1,\dots,n\}$.
Moreover, it is sufficient to take $\mathfrak{T}^{(i,j)}$ ($i\leq j$) and
$\mathfrak{T}^{(i,j,k)}$ ($i< j < k$) to obtain a generating set.
\end{thm}

The generating set of the invariant polynomial functions described in the preceding
theorem (the ones of degrees 2 and 3) is not algebraically independent. A set of
generators for their ideal of relations is given in the next theorem (second fundamental
theorem on these invariants).
\begin{thm}[see \cite{Drensky2003}, Theorem 2.3 (ii) or \cite{Drensky2006}, Theorem 3.4
(ii)]\label{th:sft}
The defining relations for the algebra of polynomial invariant functions are
\begin{subequations}\label{eq:defining_comm}
\begin{align}\label{}
 &\qquad\begin{aligned}
 &Tr([\M{i},\M{j}]\M{k}) Tr([\M{p},\M{q}]\M{r})\\
 &\quad+2\left|
 \begin{matrix}
 Tr(\M{i}\M{p}) & Tr(\M{i}\M{q}) & Tr(\M{i}\M{r})\\
 Tr(\M{j}\M{p}) & Tr(\M{j}\M{q}) & Tr(\M{j}\M{r})\\
 Tr(\M{k}\M{p}) & Tr(\M{k}\M{q}) & Tr(\M{k}\M{r})
 \end{matrix}
 \right|=0,
 \end{aligned}\\[.75em]
 &\begin{aligned}
 & Tr([\M{j},\M{k}]\M{\ell})Tr(\M{p}\M{i})
 - Tr([\M{i},\M{k}]\M{\ell})Tr(\M{p}\M{j})\\
 &+Tr([\M{i},\M{j}]\M{\ell})Tr(\M{p}\M{k})
 - Tr([\M{i},\M{j}]\M{k})Tr(\M{p}\M{\ell})=0,
 \end{aligned}
\end{align}
\end{subequations}
with $i,j,k,\ell,m,n,p,q,r\in\{1,\dots,n\}$.
\end{thm}
With the following reasoning that is adapted from
\cite{CrampePoulaindAndecyetal2020}, we now use Theorem \ref{th:fft} to extract
information about the algebra of polarized traces and the centralizer.
\begin{itemize}
\item The algebra $U(\sl_2)^{\otimes n}$ is filtered. Take the degree of all generators
$e_{ij}^{(a)}$ to be $1$, then the associated graded algebra is commutative. Recall the
$\sl_2$ action \eqref{eq:actionUsl2n}. This induces a natural action on
$gr(U(\sl_2)^{\otimes n})$. Denote the generators of the graded algebra by
$\e_{ij}^{(a)}$. The induced $\sl_2$-action is given as follows on the generators:
\begin{align}\label{eq:inducedactionUsl2n}
 e_{ij}\cdot \e_{k\ell}^{(a)}=\delta_{jk}\e_{i\ell}^{(a)}-\delta_{\ell i}\e_{kj}^{(a)}.
\end{align}
\item The algebra of polynomial functions on $\sl_2^n$ is the algebra
of polynomials  $\mathbb{C}[x_{ij}^{(a)}]$, where $x_{ij}^{(a)}$ is the linear form giving
the $(i,j)$ coordinate of the $a$\textsuperscript{th} matrix in the product
$\sl_2^n$. The simultaneous conjugation action of an element $G$ of
$SL(2)$ on a polynomial function of $\sl_2^n$ \eqref{eq:simultaneous_conjug} is given
infinitesimally by
\begin{align}\label{}
 \epsilon\sum_{k=1}^{n}f(\M{1},\dots,[\M{k},g],\dots,\M{n})
\end{align}
for $G=e^{i\epsilon g}$ with $g\in\sl_2$. Thus, on the generators of polynomials functions
$x_{ij}^{(a)}$, this infinitesimal action is
\begin{align}\label{}
 e_{ij}\cdot x_{k\ell}^{(a)}=\delta_{j\ell}x_{ki}^{(a)}-\delta_{ik}x_{j\ell}^{(a)}
\end{align}
and can be identified with the induced $\sl_2$ action on $gr(U(\sl_2)^{\otimes n})$
through
\begin{align}\label{}
\begin{aligned}
 \e_{ij}^{(a)}&\leftrightarrow x_{ji}^{(a)}.
\end{aligned}
\end{align}
\item It follows that under this identification, the invariant functions correspond to the
elements of $gr(U(\sl_2)^{\otimes n})$ in the kernel of the $\sl_2$ action, or in other
words to the image of the centralizer in $gr(U(\sl_2)^{\otimes n})$. Moreover, again under this
identification, the image in the graded algebra of the polarized trace
$T^{(a_1,\dots,a_d)}$ defined in  \eqref{eq:generatorsinvtfct} is the polynomial function
$Tr(\M{a_1}\dots \M{a_d})$.
\end{itemize}

\subsection{The defining relations of $Z_n(\sl_2)$}\label{}

Knowing that the generators of the invariant functions \eqref{eq:generatorsinvtfct}
correspond to the polarized traces \eqref{eq:polarizedtraces} (see Theorem \ref{th:fft}),
the image of the centralizer in $gr(U(\sl_2)^{\otimes n})$ is therefore generated by the
image of the polarized traces. Now consider an element of degree $N$ in the centralizer.
Up to terms of degree $N-1$, this element can be expressed as a polynomial in
$T^{(a_1,\dots,a_d)}$. The same can then be argued for each of the remaining lower degree
terms by induction, and thus any element in the centralizer can be expressed as a
polynomial in the polarized traces. Since all polarized traces belong in the centralizer,
the two algebras thus coincide. So we obtain the analogue of the first fundamental
theorem:
\begin{coro}\label{cor-FFT}
It follows from Theorem \ref{th:fft} that the polarized traces $T^{(i,j)}$, $i\leq j$ and
$T^{(i,j,k)}$, $i<j<k$ generate the centralizer $Z_n(\sl_2)$ of the diagonal action of
$U(\sl_2)$ in $U(\sl_2)^{\otimes n}$.
\end{coro}
Recall that relations \eqref{eq:defining_comm} are a set of defining relations for the
image of the centralizer in $gr(U(\sl_2)^{\otimes n})$.
We look for analogous defining relations for the centralizer in $U(\sl_2)^{\otimes n}$.
Once we find the deformations in $U(\sl_2)^{\otimes n}$ of relations
\eqref{eq:defining_comm}, we can prove their completeness using for the ideal of relations
the same sort of reasoning as before Corollary \ref{cor-FFT}. In fact the special Racah
algebra was defined such that its set of defining relations gives precisely the complete
set of relations for the centralizer.
\begin{thm}\label{thm:sft_nc}
With the following identification of the generators:
\begin{align}\label{iso}
 \p_{ab}\mapsto T^{(a,b)}\qquad\text{and}\qquad \f_{ijk}\mapsto -T^{(i,j,k)}
\end{align}
for $i,j,k$ all distinct, the centralizer $Z_n(\sl_2)$ is isomorphic to the special Racah
algebra $sR(n)$:
\begin{align}\label{}
 Z_n(\sl_2)\cong sR(n).
\end{align}
In other words, a set of defining relations of $Z_n(\sl_2)$ is given in
\eqref{eq:definingsrn}, where $\p_{ab}$ and $\f_{ijk}$ are replaced by the corresponding
polarized traces.
\end{thm}
\proof

The defining relations are verified to hold in $U(\sl_2)^{\otimes n}$ by direct
computations. Note that, due to the symmetry under renaming the indices, we only need to
make calculations in $U(\sl_2)^{\otimes n}$ for $n\leq 5$ in degrees less or equal to 6
in the generators.

Under the previous choice of degree for the generators of
$U(\sl_2)^{\otimes n}$ which was $deg(e_{ij}^{(a)})=1$, it follows that
\begin{align}\label{}
 deg(T^{(i,j)})=2,\qquad deg(T^{(i,j,k})=3.
\end{align}
The same degrees are given to the generators of the
special Racah algebra $sR(n)$: $deg(\p_{ij})=2$ and $deg(\f_{ijk})=3$. This makes it a filtered algebra, and it is straightforward
to observe that its associated graded algebra is isomorphic to the algebra of polynomial
invariants functions. Indeed, the first set of defining relations (\ref{eq:definingn}),
or equivalently \eqref{eq:srn0}--\eqref{eq:srn6},
ensures that the generators all commute in the graded algebra, and then the relations
(\ref{eq:srn7})--(\ref{eq:srn8}) and \eqref{eq:consequences} are mapped to
\eqref{eq:defining_comm}.

Therefore both algebras related by the morphism in (\ref{iso}) have the same associated graded
algebras, and so in particular have the same dimensions in each component of the
filtration (that is, in each degree). Moreover the morphism is surjective from Corollary
\ref{cor-FFT}. Consequently, an element in the kernel of the map (in other words, a relation in
$Z_n(\sl_2)$ not implied by the relations of the special Racah algebra), if non-zero,
would contradict the equality of dimensions for some degree.
\endproof

\begin{rmk}
It is quite remarkable that we only need to quotient the Racah algebra $R(n)$ by the
elements $w_{ijk}$ and $x_{ijk\ell}$ in order to recover the centralizer for any value of
$n$. Indeed, one could have expected that in order to recover the centralizer for
increasing $n$, we would need to quotient by elements of increasing degree or spanning an
increasing number of indices.  That this is not the case is quite a surprising
simplification.
\end{rmk}

\section{The Hilbert--Poincaré series of $Z_n(\sl_2)$}\label{sec:hpseries}

For more information on Hilbert--Poincaré series of graded algebras, we refer to
\cite{Stanley1978}. The Hilbert--Poincar\'e series contains useful information about a
graded, or filtered, algebra. We will illustrate this for the diagonal centralizer
$Z_n(\sl_2)$. We will provide an explicit formula for its Hilbert--Poincar\'e series, and
then use it in conjunction with the defining relations found in Theorem \ref{thm:sft_nc}
to provide bases of $Z_n(\sl_2)$ for small $n$.

\subsection{An explicit formula}

The commutative algebra $\mathbb{C}[\sl_2^n]^{inv}$ of polynomial functions on $\sl_2^n$
that are invariant under simultaneous conjugation by $SL(2)$ is a graded algebra: it is
the direct sum of the subspaces $\mathbb{C}_k[\sl_2^n]^{inv}$ of homogeneous invariant
polynomial functions of degree $k$. The Hilbert--Poincar\'e series records the dimensions
of all these subspaces:
\begin{align}\label{}
 F_n(t)=\sum_{k\geq0}\dim\bigl(\mathbb{C}_k[\sl_2^n]^{inv}\bigr)t^k.
\end{align}
The centralizer $Z_n(\sl_2)$ inherits from $U(\sl_2)^{\otimes n}$ the structure of a
filtered algebra: it is the union of the increasing sequence (in $k$) of subspaces
$Z_n(\sl_2)_{\leq k}$ of elements of degree less or equal to $k$ (the degree is in the
generators of $U(\sl_2)^{\otimes n}$). The Hilbert--Poincar\'e series of $Z_n(\sl_2)$
records the dimensions of the homogeneous subspaces of the associated graded algebra:
\begin{align}\label{}
F_n(t)=\sum_{k\geq0}\dim\bigl(Z_n(\sl_2)_{\leq k}/Z_n(\sl_2)_{<k}\bigr)t^k,
\end{align}
and thus, from the discussion in Section \ref{sec:centralizer_proof}, is the same as the
Hilbert--Poincar\'e series of the invariant polynomial functions.

Several formulas, using various approaches, have been obtained for the Hilbert--Poincar\'e
series $F_n(t)$ (see references in \cite{Drensky2006, Teranishi1986, Formanek1987}). The
formula presented below seems to be new.
\begin{prop} Let $n\geq2$ and recall that the rank $r$ is defined by $r=n-2$. The Hilbert--Poincaré
series of $Z_n(\sl_2)$  is:
\begin{align}\label{}
 F_n(t)=\frac{P_r(t)}{(1-t^2)^{3(r+1)}} ,
\end{align}
where the numerator is given by:
\begin{align}\label{}
 P_r(t)=(1+t)^r\sum_{k=0}^{2r}(-1)^ka_kt^k,\qquad\text{where}\
 \left\{\begin{array}{l}
  a_{2k}=\binom{r}{k}^2,\\[0.5em]
  a_{2k+1}=\binom{r}{k}\binom{r}{k+1}.
 \end{array}\right.
\end{align}
\end{prop}

\proof We take a detour through the graded character of $SL(2)$
on the polynomial functions on $\sl_2^n$. The character of $SL(2)$ for a finite-dimensional
representation is seen as a Laurent polynomial in $x$, given by the trace of the action of
the element $\text{Diag}(x,x^{-1})$. For example, for the fundamental representation of
$SL(2)$, it is $x+x^{-1}$. For the irreducible representation of dimension $d+1$, which is
the $d$-symmetrized power of the fundamental representation, the character is thus
$\frac{x^{d+1}-x^{-d-1}}{x-x^{-1}}$.  Now, it is easy to check that, if we have the
character $\chi(x)$ of an arbitrary finite-dimensional representation of $SL(2)$, then the
formula:
\begin{align}\label{}
 \bigl[(1-x^2)\chi(x)\bigr]_0,
\end{align}
where $[\,\cdot\,]_0$ means taking the constant term of a Laurent polynomial, gives the
multiplicity of the trivial representation.

After these classical preliminaries, note that the character of the adjoint representation
of $SL(2)$ on $\sl_2$ is $1+x^2+x^{-2}$. On the polynomial function on $\sl_2$, the action
of $SL(2)$ preserves the grading, and we record the character of the representation on its
graded components as a formal power series in $t$ (also called, the graded character). For
each degree, the representation is a symmetrized power of the adjoint representation, so
we find that the graded character is:
\begin{align}\label{}
 \frac{1}{(1-t)(1-tx^2)(1-tx^{-2})}.
\end{align}
Equivalently, this is the graded character for the adjoint action on $U(\sl_2)$. On the
polynomial functions $\sl_2^n$ (or equivalently, on $U(sl_2)^{\otimes
n}$), the graded character is thus:
\begin{align}\label{}
 \frac{1}{(1-t)^n(1-t x^2)^n(1-t x^{-2})^n}.
\end{align}
Now, in each degree, we look for the dimension of the invariant subspace for the action of
$SL(2)$. In other words, we look for the multiplicity of the trivial representation. By
what we have recalled above, we obtain that the Hilbert--Poincar\'e series of $Z_n(\sl_2)$
is:
\begin{align}\label{}
 F_n(t)=\left[\frac{(1-x^2)}{(1-t)^n(1-t x^2)^n(1-t x^{-2})^n}\right]_0.
\end{align}
Using the expansion $(1-z)^{-n}=\sum_{k\geq 0}\binom{k+n-1}{k}z^k$ and straightforward
manipulations, we obtain:
\begin{equation}\label{eq-inter1}
 F_n(t)=\frac{1}{(1-t)^n}\sum_{k\geq0}(-1)^k\tilde{a}_kt^k,\qquad\text{where}\
 \left\{\begin{array}{l}
  \tilde{a}_{2k}=\binom{n+k-1}{k}^2 ,\\[0.5em]
  \tilde{a}_{2k+1}=\binom{n+k-1}{k}\binom{n+k}{k+1} .
 \end{array}\right.
\end{equation}
Thus the statement of the proposition reduces to the following equality of formal power
series:
\begin{equation}\label{eq-series}
 \frac{1}{(1-t)^{r+2}}\sum_{k\geq0}(-1)^k\tilde{a}_kt^k
 =\frac{(1+t)^r}{(1-t^2)^{3(r+1)}}\sum_{k\geq0}(-1)^ka_kt^k.
\end{equation}
To prove this, we multiply the $\tilde{a}$-series by $(1+t)$ and the $a$-series by
$(1-t)$, and after an application of Pascal rule for binomials, we reach the equivalent
formula:
\begin{align}\label{}
 \sum_{k\geq0}(-1)^k\tilde{a}'_kt^k=\frac{1}{(1-t^2)^{2(r+1)}}\sum_{k\geq0}(-1)^ka'_kt^k,
\end{align}
where now we have
$\left\{\begin{array}{l}
 \tilde{a}'_{2k}=\binom{r+k+1}{k}\binom{r+k}{k},\\[0.5em]
 \tilde{a}'_{2k+1}=\binom{r+k+1}{k}\binom{r+k+1}{k+1},
 \end{array}\right.$
and\quad
$\left\{\begin{array}{l}
 a'_{2k}=\binom{r}{k}\binom{r+1}{k} ,\\[0.5em]
 a'_{2k+1}=\binom{r}{k}\binom{r+1}{k+1} .
\end{array}\right.$\\[.5em]
This last formula is verified by writing the expansion of the right hand side and checking
the equality of the coefficients, making use of the following identity for binomial
coefficients \cite{Gould1972}:
\begin{align}\label{}
 \sum_{i}\binom{i+a+b}{i}\binom{b}{k-i}\binom{a}{k'-i}=\binom{k'+b}{k}\binom{k+a}{k'}.
\end{align}
This identity is valid for any $a,b,k,k'$ and we use it for $a=r+1$, $b=r$ and
$k'\in\{k,k+1\}$.
\endproof


The exponent $3(r+1)$ appearing in the denominator of $F_n(t)$ is the Krull, or
Gelfand--Kirillov, dimension of the algebra of invariant polynomial functions (see
\cite{Kirillov1967},\cite{Formanek1987}). Here it means that there is a set of $3(r+1)$
algebraically independent homogeneous elements (a system of parameters)
$\theta_1,\dots,\theta_{3(r+1)}$ such that the algebra is a free module of finite
dimension over the polynomial subalgebra $\mathbb{C}[\theta_1,\dots,\theta_{3(r+1)}]$. The
freeness follows from the property called Cohen--Macaulay, which is ensured here by the
Hochster--Roberts theorem from general invariant theory \cite{HochsterRoberts1974}. The
form $F_n(t)$ above with the positivity of the numerator $P_r(t)$ (see below) suggests that it might
be possible that a system of parameters consists of $3(r+1)$ elements of degrees $2$. If it
were to be the case, then $P_r(1)$ would be the dimension of the algebra over
$\mathbb{C}[\theta_1,\dots,\theta_{3(r+1)}]$. Moreover, the different monomials in
$P_r(t)$ would indicate in which degrees the elements of a basis over
$\mathbb{C}[\theta_1,\dots,\theta_{3(r+1)}]$ would have to be found.

Finally, the palindromic property of the numerator $P_r(t)$ in the formula above shows
that the Hilbert--Poincar\'e series satisfies the functional equation:
\begin{align}\label{}
 F_n(t^{-1})=(-1)^{(n-1)}t^{3n}F_n(t).
\end{align}
This is well-known and related to a property, called being Gorenstein, for the algebra of
invariant polynomial functions (see \cite{Drensky2006} and references therein).

\begin{rmk}
More generally, the Hilbert--Poincar\'e series of $Z_n(\sl_2)$ can be defined as a power
series in $t_1,\dots,t_n$ if we consider the gradation by the multidegree of
$U(\sl_2)^{\otimes n}$. A slight generalization of the part of the proof up to Formula
(\ref{eq-inter1}) gives the multigraded version of this formula:
\begin{align}\label{}
 F_n(t_1,\dots,t_n)=\frac{1}{(1-t_1)\dots(1-t_n)}
 \Bigl(\sum_{\substack{\mu\models k\\ \nu\models k}}
  t_1^{\mu_1+\nu_1}\dots t_n^{\mu_n+\nu_n}
 -\sum_{\substack{\mu\models k\\ \nu\models k-1}}
  t_1^{\mu_1+\nu_1}\dots t_n^{\mu_n+\nu_n}\Bigr),
 \end{align}
where $\mu\models k$ means that $\mu=(\mu_1,\dots,\mu_n)\in\mathbb{Z}_{\geq0}^n$ such that
$\mu_1+\dots+\mu_n=k$. To recover relation (\ref{eq-inter1}) from this, take
$t_1=\dots=t_n=t$ and use that the number of $\mu\models k$ is $\binom{k+n-1}{k}$.
\end{rmk}

\subsection{Some related combinatorics}

We have obtained an expression for the Hilbert--Poincar\'e series of $Z_n(\sl_2)$ of the
form:
\begin{align}\label{}
 F_n(t)=\frac{(1+t)^rQ_r(t)}{(1-t^2)^{3(r+1)}} ,\qquad
 \text{where $Q_r(t)=\sum_{k=0}^{2r}(-1)^ka_kt^k$},
\end{align}
and the coefficients $a_k$ are given in the proposition above. It is perhaps not so
surprising that the coefficients of the various polynomials involved show some connections
with well-studied combinatorial objects of ``Catalan'' flavor.

\paragraph{The polynomial $Q_r(t)$.} The coefficient $a_k$ in the polynomial $Q_r(t)$
counts the number of symmetric Dyck paths of semi-length $2r+1$ with $k+1$ peaks (see
A088855 in \cite{OEIS}). Their expression with binomial coefficients corresponds to
choosing a certain number of peaks and troughs in the first $r$ steps of the paths.

In fact, the polynomial $Q_r(t)$ is a $t$-deformation of the well-known Catalan number,
that is, the value of $Q_r(t)$ at $t=1$ is the $r$-th Catalan number:
\begin{align}\label{}
 Q_r(1)=c_r=\binom{2r}{r}-\binom{2r}{r+1}.
\end{align}
Indeed it is not difficult to give a combinatorial proof that the alternating sum of the
$a_k$'s is equal to the Catalan number $c_r$ (the number of Dyck paths of length $2r$). We
can also see it as follows. The Catalan number is equal to the constant term of a Laurent
polynomial in $x$:
\begin{align}\label{}
 c_r=\bigl[(1-x^2)(x+x^{-1})^{2r}\bigr]_0
    =\bigl[(1-x^2)(1+x^{2})^{r}(1+x^{-2})^{r}\bigr]_0.
\end{align}
Note that we keep the variable $x^2$ to stay coherent with the notation used during the
proof above. In fact, the first equality expresses that the Catalan number $c_r$ is the
multiplicity of the trivial representation in $V^{r}\otimes (V^{\star})^r$. The
$t$-deformation is now immediate from this formula for $c_r$. Indeed, it follows from its
explicit expression that the polynomial $Q_r(t)$ is given by:
\begin{align}\label{}
 Q_r(t)=\bigl[(1-x^2)(1+tx^{2})^{r}(1+tx^{-2})^{r}\bigr]_0.
\end{align}
In this sense, the polynomial $Q_r(t)$ is a natural $t$-deformation of the $r$-th Catalan
number.

\paragraph{The numerator $P_r(t)$.} The numerator of the Hilbert--Poincar\'e series of
$Z_n(\sl_2)$ is $P_r(t)=(1+t)^rQ_r(t)$. We will show explicitly that its coefficients are
positive.

First, from what we have said above about $Q_r(t)$, it follows that $P_r(t)$ is a
$t$-deformation of the number $2^rc_r$, that is, its value at $t=1$ is $P_r(1)=2^rc_r$.
This number counts several classes of combinatorial objects (obtained from objects counted
by the Catalan number, see A151374 in \cite{OEIS}). The $t$-deformation giving $P_r(t)$
can be expressed similarly as before as:
\begin{align}\label{}
 P_r(t)=\bigl[(1-x^2)(1+t)^r(1+tx^{2})^{r}(1+tx^{-2})^{r}\bigr]_0.
\end{align}
Now regrouping the terms with an $r$-th power, this gives the following expression:
\begin{align}\label{form-P}
\begin{aligned}
 P_r(t) & =  \bigl[(1-x^2)\bigl(1+t^3+(t+t^2)(1+x^2+x^{-2})\bigr)^{r}\bigr]_0\\
 & = \displaystyle\sum_{k=0}^rR_k\binom{r}{k}(1+t^3)^{r-k}(t+t^2)^k,
\end{aligned}
\end{align}
where the positive integer $R_k$ is the Riordan number, one of the closest relative of the
Catalan number, which also admits many combinatorial interpretations (see A005043 in
\cite{OEIS}). They are given by either one of the following equalities:
\begin{align}\label{}
 R_n=\bigl[(1-x^2)(1+x^2+x^{-2})^n\bigr]_0=\sum_{i=0}^n(-1)^{n-i}\binom{n}{i}c_i.
\end{align}
We actually used the first one in the above calculation, while the second one shows that
the Catalan sequence  is the binomial transform of the Riordan sequence, and thus allows
to recover that $P_r(1)=2^rc_r$.

The formula (\ref{form-P}) for $P_r(t)$ has the advantage to show explicitly that it has
positive coefficients. So $P_r(t)$ is a $t$-deformation with positive coefficients of
$2^rc_r$ and therefore should be given by an interesting statistics on a certain set of
$2^rc_r$ objects.

\subsection{PBW basis of $Z_n(\sl_2)$ for small $n$}

With the help of the Poincar\'e--Hilbert series obtained and discussed above,
we here determine the PBW bases of $Z_n(\sl_2)$ for $n=2,3,4$.
Using Theorem \ref{thm:sft_nc}, we identify systematically $Z_n(\sl_2)$ with the special
Racah algebra, and use the generators and relations in Section \ref{sec:sRacah}. Note that
from the discussion in Section \ref{sec-inv}, in order to show that a subset of
$Z_n(\sl_2)$ is a spanning set, it is enough to show that the images in the graded algebra
of invariant polynomial functions is a spanning set.

\paragraph{The case $n=2$.} The Hilbert--Poincar\'e series of $Z_2(\sl_2)$ is:
\begin{equation}\label{Hn-2}
F_2(t)=\frac{1}{(1-t^2)^{3}} .
\end{equation}
This allows easily to recover that the algebra of invariant polynomial functions, and
$Z_2(\sl_2)$, is a commutative polynomial algebra. A basis of  $Z_2(\sl_2)$ is
$\p_{11}^{a}\p_{12}^{b}\p_{22}^{c}$ with $a,b,c\in\mathbb{Z}_{\geq0}$. Indeed this set is
obviously a spanning set, and moreover spans a subspace whose dimensions in each degree
are given precisely by the series (\ref{Hn-2}). So this set is also linearly independent,
and is a basis.

\paragraph{The case $n=3$.} The Hilbert--Poincar\'e series of $Z_3(\sl_2)$ is:
\begin{align}\label{}
 F_3(t)=\frac{1+t^3}{(1-t^2)^{6}}.
\end{align}
This allows to show that the following set is a basis of $Z_3(\sl_2)$:
\begin{equation}\label{basis-H3}
 \{ \p_{11}^{a}\p_{12}^{b}\p_{13}^c\p_{22}^{d}\p_{23}^{e}\p_{33}^{f}\f_{123}^{\varepsilon} \},
 \qquad\text{where $a,b,c,d,e,f\in\mathbb{Z}_{\geq0}$ and $\varepsilon\in\{0,1\}$.}
\end{equation}
Indeed, such a set is a spanning set since ${\f_{123}}^2$ can be expressed in terms of the
$\P$'s. Comparing with $F_3(t)$, we see directly that this set spans a subspace of the
correct dimension in each degree. So this is a basis.

\paragraph{The case $n=4$.} The Hilbert--Poincar\'e series of $Z_4(\sl_2)$ is:
\begin{align}\label{}
 F_4(t)=\frac{1+t^2+4t^3+t^4+t^6}{(1-t^2)^9}.
\end{align}
The four sets, for $a,b,c,d,e,f,g,h,i \in \mathbb{Z}_\geq 0$,
\begin{subequations}\label{eq:pbw_sr4}
\begin{align}\label{basis-H4}
& \{ \f_{123}\p_{11}^{a}\p_{12}^{b}\p_{13}^c\p_{14}^d
\p_{22}^{e}\p_{23}^{f}\p_{24}^{g}\p_{33}^{h}\p_{34}^{i} \},   \\
& \{ \f_{124}\p_{11}^{a}\p_{12}^{b}\p_{13}^c\p_{14}^d \p_{22}^{e}\p_{23}^{f}\p_{24}^{g}
\p_{34}^{h}\p_{44}^{i} \},  \\
& \{ \f_{134}\p_{11}^{a}\p_{12}^{b}\p_{13}^c\p_{14}^d
\p_{23}^{e}\p_{24}^{f}\p_{33}^{g}\p_{34}^{h}\p_{44}^{i} \}, \\
&  \{ \f_{234} \p_{12}^{a}\p_{13}^b\p_{14}^c
\p_{22}^{d}\p_{23}^{e}\p_{24}^{f}\p_{33}^{g}\p_{34}^{h}\p_{44}^{i} \},
\end{align}
and, for $a,b,c,d,e,f,g,h,i,j \in \mathbb{Z}_\geq 0$ and $aehj=0$,
\begin{align}
& \{ \p_{11}^{a}\p_{12}^{b}\p_{13}^c\p_{14}^d
\p_{22}^{e}\p_{23}^{f}\p_{24}^{g}\p_{33}^{h}\p_{34}^{i}\p_{44}^{j} \}
\end{align}
\end{subequations}
form a basis of $Z_4(\sl_2)$.
To understand this, rewrite its Hilbert--Poincaré series as
\begin{equation}
 F_4(t)= \frac{4t^3}{(1-t^2)^9} + \frac{1-t^8}{(1-t^2)^{10}}.
\end{equation}
The first term corresponds to the first 4 sets and the second term corresponds to the fifth ones.
These sets are spanning sets since $\f_{ijk} \f_{mnp} $ can be expressed in terms linear in $\f$
by using the relations $w_{ijk}=0$ and $x_{ijk\ell}=0$ and $\f_{ijk} \p_{\ell\ell} $ ($i,j,k,\ell$ pairwise distinct)
can be expressed in terms of elements of the sets \eqref{eq:pbw_sr4} by using \eqref{eq:r4d5}.
The condition $aehj=0$ for the fifth set comes from the following fact.
Let us define the following $2\times2$ and $4\times4$ matrices
\begin{align}\label{}
 \p_{ij}^{\,ab}=
 \begin{pmatrix}
  \p_{ia} & \p_{ib} \\
  \p_{ja} & \p_{jb}
 \end{pmatrix},\qquad
 \p_{ijk\ell}^{\,abcd}=
 \begin{pmatrix}
  \p_{ia} & \p_{ib} & \p_{ic} & \p_{id}\\
  \p_{ja} & \p_{jb} & \p_{jc} & \p_{jd}\\
  \p_{ka} & \p_{kb} & \p_{kc} & \p_{kd}\\
  \p_{\ell a} & \p_{\ell b} & \p_{\ell c} & \p_{\ell d}
 \end{pmatrix}.
\end{align}
Recall the definition of the symmetrized determinant \eqref{eq:sym_det}.
In  $Z_4(\sl_2)$, the following relation of degree $8$ is satisfied by the generators:
\begin{align}\label{}
 \det(\p_{1234}^{1234})&=
 -\tfrac13\left( \det(\p_{123}^{124})-\det(\p_{123}^{134})+\det(\p_{123}^{234})
           +\det(\p_{124}^{134})-\det(\p_{124}^{234})+\det(\p_{134}^{234})\right)\nonumber \\
 &\quad+\tfrac23\left( \p_{12}\det(\p_{34}^{34})+\p_{13}\det(\p_{24}^{24})
                 +\p_{14}\det(\p_{23}^{23})\right. \nonumber \\
 &\hspace{3em}  \left. +\p_{23}\det(\p_{14}^{14})+\p_{24}\det(\p_{13}^{13})
                 +\p_{34}\det(\p_{12}^{12})\right).
\end{align}
Note that the above relation is not a new relation (it is implied by the defining
relations of $sR(4)$ given in \eqref{eq:definingsrn}--\eqref{eq:consequences}) and permits to express $\p_{11}\p_{22}\p_{33}\p_{44}$ in terms of the elements of the sets \eqref{eq:pbw_sr4}.


\section{Conclusion}\label{sec:concl}

Classical results about the invariant theory of the polynomials on $\mathfrak{sl}_2^n$
has allowed to provide a description in terms of generators and relations of the diagonal
centralizer of $\mathfrak{sl}_2$. A precise connection with the higher rank Racah algebra was given.
Various questions arise and pave the way to
different generalizations.

We would like to emphasize that the natural numbers appearing in the numerator of the
Hilbert--Poincaré form very well-known series of integers that have numerous interpretations
and appear already in the study of the representation theory of $\mathfrak{sl}_2$.  This
suggests that there should be a further understanding of these numbers.

The classification of the finite irreducible representations of the rank 1 Racah
algebra has been done in \cite{Huang2021}. Following this, it should be possible to
study the finite-dimensional representations of the (special) higher rank Racah algebra
$R(n)$ (resp. $sR(n)$).  These representations must be closely related to the
operators associated to the $(n-2)$-variable Racah polynomials. Indeed, the difference and
recurrence operators characterizing the univariate Racah polynomials satisfy the
relations of $R(3)$.  The study of the generalization to $(n-2)$-variable polynomials has
been initiated in \cite{DeBieIlievetal2020, DeBievandeVijver2020} and it would be
interesting to verify if the operators used in this case realize $R(n)$ or $sR(n)$.  The
Racah algebra appears also as the symmetry algebra of some superintegrable models
\cite{KalninsKressetal2005, Post2011a, KalninsMilleretal2013, GenestVinetetal2014a,
GenestVinetetal2014}. We trust that the observations and theorems of the present paper will lead to
a deeper understanding of this symmetry.

We focused here on the diagonal centralizer of $\mathfrak{sl}_2$ in the $n$-fold tensor
product of $\mathfrak{sl}_2$.  Other cases where $\sl_2$ is replaced by other algebras are
also known.  For example the diagonal centralizer of the oscillator algebra has been
studied in \cite{CrampevandeVijveretal2020}, the diagonal centralizer of the super Lie
algebra $\mathfrak{osp}(1|2)$ is known to be related to the Bannai--Ito algebra \cite{Genestetal2012}
and the centralizer of $\mathfrak{sl}_3$ in its twofold tensor product has been introduced
in \cite{CrampePoulaindAndecyetal2020}.  An important generalization concerns the quantum
group $U_q(\mathfrak{sl}_2)$.  Its diagonal centralizer in the $3$-fold tensor product has
been examined \cite{Zhedanov1991} and is associated to the Askey--Wilson algebra (see e.g.
\cite{CrampeFrappatetal2021} for a review).  Many attempts
\cite{PostWalter2017,DeBieDeClercqetal2020, DeClercq2019} to generalize this result to
$n$-fold tensor products have yielded relations of this centralizer but certainly did not
give all the defining relations.  Looking ahead, we are planning to provide a complete set of
defining relations, by using some deformation of the defining relations of the special
Racah algebra given in this paper.

We have studied the centralizer $Z_n(\sl_2)$ at an algebraic level. It is however equally important to study this centralizer
when each factor in the $n$-fold tensor product is in a finite-dimensional irreducible representation of
$\mathfrak{sl}_2$. In the case
$n=3$, a conjecture stating that the centralizer is a quotient of the Racah algebra
$R(3)$ was given in \cite{CrampePoulaindAndecyetal2019} (see \cite{CrampeVinetetal2020a} for the $q$-deformed case
and \cite{CrampeFrappatetal2019} for
$\mathfrak{osp}(1|2)$). This quotient associates the Racah
algebra with well-known algebras such as the Temperley--Lieb or Brauer algebras.
The generalization of these results to the case of the
$n$-fold tensor product is desirable and the results obtained in the present paper offer a nice
starting point.  As another follow-up, we plan on finding the explicit quotient that provides a
description of these centralizers in representations and to compare them to the recent results reported in
\cite{CrampePoulaindAndecy2020,FloresPeltola2020}.

\subsection*{Acknowledgments}

N. Crampé and L. Poulain d'Andecy are partially supported by Agence Nationale de la
Recherche Projet AHA ANR-18-CE40-0001.
J. Gaboriaud held an Alexander-Graham-Bell scholarship from the
Natural Sciences and Engineering Research Council of Canada (NSERC) and received scholarships from the ISM and the Université de Montréal.
The research of L. Vinet is funded in part by a Discovery Grant from NSERC.
\\

\printbibliography

\end{document}